\newcounter{version}
\documentstyle[12pt]{amsart}
\numberwithin{equation}{section}

\ifnum\theversion>0
  
\fi

\newtheorem{thm}{Theorem}[section]

\newtheorem{lem}[thm]{Lemma}

\theoremstyle{definition}
\newtheorem{defn}[thm]{Definition}

\theoremstyle{remark}
\newtheorem{rem}{Remark}

\newcommand{\D}{\Delta}
\newcommand{\Un}{\mathop{\text{U}}\nolimits}

\newcommand{\Proj}{{\bf P}}
\newcommand{\R}{{\Bbb R}}
\newcommand{\im}{\sqrt{-1}}

\newcommand{\Spin}{\mathop{\text{SPIN}}\nolimits}

\newcommand{\SO}{\mathop{\text{SO}}\nolimits}

\begin{document}

\newcommand{\Pin}{\mathop{\rm PIN}\nolimits}
\newcommand{\spin}{\mathop{\frak{spin}\,}\nolimits}
\newcommand{\un}{\mathop{{\frak u}\,}\nolimits}
\newcommand{\so}{\mathop{\frak{so}}\nolimits}
\newcommand{\C}{{\Bbb C}}
\newcommand{\Z}{{\Bbb Z}}
\newcommand{\SS}[1]{s^{\overline {#1}}}
\newcommand{\ET}[2]{\eta^{\overline {#1} #2}}
\newcommand{\AAA}{{\cal Q}_m}
\newcommand{\FF}{{\cal F}}
\newcommand{\partdel}[2]{\frac{\partial#1}{\partial#2}}
\newcommand{\liedela}[1]{L_{\partdel{}{\xi^#1}}}
\newcommand{\liedelb}[1]{L_{\overline{\partdel{}{\xi^#1}}}}
\newcommand{\Zpm}{Z^{\pm}}
\newcommand{\Zp}{Z^{+}}
\newcommand{\Zm}{Z^{-}}
\newcommand{\vfld}{{\cal F}}
\newcommand{\partdelx}[1]{\partdel{}{x_{#1}}}
\newcommand{\End}{\text{End}}
\newcommand{\symindex}{I_1,\ldots,I_m}
\newcommand{\sss}{\overline{s^{\symindex}}}
\newcommand{\tttt}{\overline{t^{\symindex}}}
\newcommand{\idw}[1]{i(\partdel{}{\overline{w_{#1}}})}
\newcommand{\idwij}{\idw{ij}}
\newcommand{\idwkl}{\idw{kl}}
\newcommand{\transpose}[1]{{}^t#1}
\newcommand{\symtheta}{\theta_{\symindex}}
\newcommand{\twoform}[2]{e^{#1}\wedge e^{#2}}
\newcommand{\eee}{{\cal E}}
\newcommand{\Dp}{\D\!^+}
\newcommand{\Dm}{\D\!^-}
\newcommand{\Dpm}{\D\!^\pm}
\newcommand{\Dmp}{\D\!^\mp}
\newcommand{\emp}{\emptyset}
\newcommand{\dbar}{\overline{\partial}}
\newcommand{\CO}{\mathop{\text{CO}}\nolimits}
\newcommand{\CE}{\mathop{\text{CE}}\nolimits}
\newcommand{\PPP}{{\cal P}}

\title[The inverse Penrose transform]
{The inverse Penrose transform
	on Riemannian twistor spaces}
\author{Yoshinari Inoue}
\address{Department of Mathematics\\
	Faculty of Science\\ Kyoto University\\
	Kyoto 606-01, Japan}
\email{inoue@@kusm.kyoto-u.ac.jp}

\maketitle

\begin{abstract}
With respect to the Dirac operator and the conformally invariant
Laplacian, an explicit description of the inverse Penrose transform
on Riemannian twistor spaces is given.
A Dolbeault representative of cohomology on the twistor space
is constructed from a solution of the field equation on the base manifold.
\end{abstract}

\bigbreak
\noindent
{\bf Introduction}
\bigbreak

The Penrose transform is a method to give solutions
of the Dirac equation and the conformally invariant Laplacian.
It is done by relating 
  solutions of the field equations
to cohomology
with values in a certain holomorphic line bundle
over the twistor space of the manifold.

The Penrose transform on four dimensional half-conformally
flat manifolds was
studied by Hitchin in~[H].
Murray generalized it to higher dimensional
conformally flat manifolds
in [M]. 

The correspondence between cohomology
and the space of solutions of the field equation was
proved to be one-to-one.
But the sufficiency part of the proof, that is,
to construct a cohomology class from
a solution of the field equation,
was proved indirectly in both papers.
In the four dimensional case, an explicit formula
for the inverse Penrose transform was given by
Woodhouse in [W].

In this paper, we shall give an explicit description of
the inverse Penrose transform for Riemannian manifolds by
constructing a Dolbeault representative
of the corresponding cohomology
(Definition~\ref{def4} in the four dimensional case
and  Definition~\ref{InvPenTr}  in the higher dimensional cases).
Definition~\ref{def4} is equivalent to the formula
given by Woodhouse.

Let $M$ be a $2n$-dimensional spin manifold
and let $V$ be a Hermitian vector bundle on $M$ with a connection.
We assume the conditions on the metric of $M$
and the curvature of $V$ which enable us to
perform the Penrose transform.
Then the twistor space $Z^\pm(M)$ of $M$
is a complex manifold, and the hyperplane bundle $H$
and the pull-back bundle $p^{-1}V$ on $Z^\pm(M)$
are holomorphic bundles.
Let $\Dpm(M)$ be the spin bundle on $M$.
Then the twisted differential form on
 $Z^\pm(M)$ which represents
the cohomology class corresponding to
$\phi\in\Gamma(M, S^m\Dpm(M)\otimes V)$  is written as
$$
	\AAA(\phi) = (n+m-2)! F^{(n+m-2)}(D)j(\phi),
 $$
where $D$ is a differential operator and $F^{(n+m-2)}$ is
the $(n+m-2)$-th derivative of the power series
$F(x)=\sum_{k=0}^\infty x^k/ (k!)^2$.
The lifting $j(\phi)$ is written as
$$
        j : \Gamma(M,S^m\D^{\pm}(M)\otimes V) \rightarrow
            \Gamma(\Zpm(M),
                \Lambda_V^{0,\frac{1}{2}n(n-1)}Z^\pm(M)\otimes H^{-2n+2-m}
		  \otimes p^{-1}V),
 $$
where $\Lambda_V^{0,\frac{1}{2}n(n-1)}Z^\pm(M)$
is a line subbundle of $\Lambda^{0,\frac{1}{2}n(n-1)}Z^\pm(M)$.
The definitions of $D$ and 
 $\Lambda_V^{0,\frac{1}{2}n(n-1)}Z^\pm(M)$
in the four dimensional case
 are slightly different from
those in the higher
dimensional cases.

In the four dimensional case,
 $\Lambda_V^{0,1}Z^\pm(M)$
is the space of vertical forms with respect to the
Levi-Civita connection on $M$, which is defined globally.
The operator $D$ is also defined as a global operator.

In the higher dimensional cases, the assumption of the metric of $M$
means that there are local conformally flat coordinates.
The line subbundle $\Lambda_V^{0,\frac{1}{2}n(n-1)}Z^\pm(M)$
is the space of vertical forms corresponding to the trivialization
with respect to these coordinates.
The operator $D$ is also defined with respect to them.
Although $j$ and $D$ are defined to be local operators,
it is shown that
the constructed inverse Penrose transform is defined
independently of the particular conformally flat
coordinates which are chosen.

In both cases, it is shown that $D^{n+1}$ vanishes.
Hence the construction
is in fact a finite sum.

In the  proof of the vanishing of $\dbar\AAA(\phi)$ and
the independence of the construction of $\AAA$ with respect to the
coordinates chosen
in the higher dimensional cases,
the defining equations of $Z^\pm$ as a subvariety of $\Proj(\Dpm)$
play an important role
(see Lemma \ref{spin_rel} below).
They were given in [I1]
in the course of my definition of Riemannian twistor spaces.
Although they are trivial in four and six
dimensional cases, they are still useful by
our extension of notation of multi-indices.


In \S1, we review the theory of the Penrose transform
on Riemannian twistor spaces.
In \S2 and \S3, we give the inverse Penrose transform
on $2n$-dimensional manifolds with $n\ge3$.
The local construction of the inverse Penrose transform
is given in \S2,
 and 
the independence of the construction with respect to
the coordinates which are chosen
is proved in \S3.
In \S4, we deal with the four dimensional case,
in which a non-conformally flat manifold may admit the
integrable twistor space.

I would like to thank the reviewer for his
suggestion of the work of Woodhouse.

\section{The Penrose transform}

Let us review the Penrose transform for the Dirac operator and
the conformally invariant Laplacian
on even dimensional spin manifolds (see [H]: the four dimensional case,
 and [M]: the higher dimensional cases).

Let $M$ be a $2n$-dimensional
spin manifold, and let $V$ be a Hermitian vector bundle on $M$
with a connection.
Let $\Dp(M)$ (resp. $\Dm(M)$) be the positive (resp. negative) spin bundle. 
The field equations that we consider are
the Dirac operator
$$
        d_m: \Gamma(M,S^m\Dpm(M)\otimes V) \rightarrow
                \Gamma(M,\Dmp(M)\otimes S^{m-1}\Dpm(M)\otimes V),
 $$
for $m>0$ and the conformally invariant Laplacian
\begin{align*}
        d_0: \Gamma(V) &\rightarrow \Gamma(V)\\
               \phi &\mapsto \nabla^*\nabla(\phi) +
		 \frac{n-1}{2(2n-1)}r\phi,
\end{align*}
where $r$ is the scalar curvature of $M$.
The differential operator $d_m$
is conformally invariant
with conformal weight $n-1 +\frac{m}{2}$.

Assume $n\ge2$.
Let $Z^{+}$ be the parameter space of complex structures  of
a $2n$-dimensional real vector space
compatible with a certain metric and a certain orientation.
If we consider the opposite orientation, we have a similar manifold $Z^{-}$.
Let $\SO(M)$ be the oriented orthonormal frame bundle of $M$.
Riemannian twistor spaces are defined as
$$
	Z^{\pm}(M) = \SO(M)\times_{\SO(2n)}Z^{\pm},
 $$
which have natural almost complex structures.
We assume that $M$ has an integrable twistor space.
Let $p$ be the projection
$p: Z^\pm(M) \rightarrow M $.
Assume that the pull-back of the curvature form of $V$ by $p$
is a $\End(p^{-1}V)$-valued $(1,1)$-form.
Then the pull-back bundle $p^{-1}V$ can be
naturally considered to be a holomorphic vector bundle.

There is a holomorphic line bundle  $H$ on $Z^\pm(M)$,
which is called the hyperplane bundle.
Let us define the Penrose transform:
$$
	\PPP_m : H^{\frac{1}{2}n(n-1)}(Z^{\pm}(M),H^{-2n+2-m}\otimes p^{-1}V)
		\rightarrow \Gamma(M,S^m\Dpm(M)\otimes V).
 $$
By restricting  cohomology classes to
each fiber, we have a map
$$
	H^{\frac{1}{2}n(n-1)}(Z^{\pm}(M),H^{-2n+2-m}\otimes p^{-1}V)
	\rightarrow \bigcup_{x\in M}
	  H^{\frac{1}{2}n(n-1)}(Z^{\pm}_x,H^{-2n+2-m})\otimes V_x.
 $$
This induces $\PPP_m$,
since  $H^{\frac{1}{2}n(n-1)}(Z^{\pm}_x,H^{-2n+2-m})$
is equivalent to $S^m\D^{\pm}_{\,x}$ as a representation space
of $\Spin(T_xM)$ by the theorem of Bott-Borel-Weil-Kostant.

\begin{thm} [{[H] Theorem (3.1),  [M] Theorem 32}]
The map $\PPP_m$ induces
 an isomorphism onto the space of the solutions of $d_m\phi=0$.
\end{thm}

\section{ The Local construction of the inverse Penrose transform}
\label{const}

In this section, we deal with the local construction
of the inverse Penrose transform when the base manifold
is $2n$-dimensional with $n\ge3$.
We assume  that $M$ is an open subset of 
$\R^{2n}$ with the standard metric
and the vector bundle $V$ is trivial.
Hence we can safely  omit $V$
by considering it as a trivial line bundle.

By Dolbeault's theorem, we have a representation of the cohomology group:
$$
	H^{\frac{1}{2}n(n-1)}(Z^{\pm}(M),H^{-2n+2-m})
	  =
	\frac{ \ker
	  \dbar|_{
	  \Gamma(\Lambda^{0,\frac{1}{2}n(n-1)}
	     \Zpm(M)\otimes H^{-2n+2-m})}
	   }
	   {\dbar\Gamma(\Lambda^{0,\frac{1}{2}n(n-1)-1}
             \Zpm(M)\otimes H^{-2n+2-m})}.
 $$
We  construct a  $\dbar$-closed form in
$\Gamma(\Zpm(M), \Lambda^{0,\frac{1}{2}n(n-1)}
             \Zpm(M)\otimes H^{-2n+2-m})$
from a solution of $d_m\phi=0$.

We begin by defining a differential operator $D$
acting on $\Gamma(\Zpm(M), \Lambda^{\!*}
\Zpm(M)\otimes H^{-2n+2-m})$.

The action of $\SO(2n)$ on $\Zpm$ induces a linear map
$$
	\cal{F}: \so(2n) \rightarrow \Gamma(\Zpm,\Theta),
 $$
where $\Theta$ is the holomorphic tangent bundle on $\Zpm$, which
is considered to be the $(1,0)$-part of the complexified tangent
bundle $T\Zpm\otimes\C$. With respect to the Lie algebra structures, we have
\begin{equation}\label{lie_alg_str}
	{\cal F}([a,b]) = - [{\cal F}(a),{\cal F}(b)].
\end{equation}
Let
$$
  \begin{aligned}
	e_a&=\partdel{}{x_a},\\
	e^a&=dx_a,
  \end{aligned}
	\qquad a=1,\ldots,2n
 $$
be the standard frames of $TM$ and $T^*M$, respectively.
Let $(E^a_b)_{1\le a,b\le 2n}$ be the frame of $\End(TM)$ defined as
$E^a_be_c= \delta^a_ce_b$
and put $F_{ab}= E^a_b - E^b_a$.
By considering it as an element of $\so(2n)$,
we define a vector field $\vfld_{ab}$ on $M\times Z^{\pm}=Z^{\pm}(M)$ to be
\begin{equation*}
	\vfld_{ab} = \vfld(F_{ab}).
\end{equation*}
Then we define a first-order differential operator $D$ acting on
 $\Gamma(\Zpm(M), \Lambda^{\!*} \Zpm(M)\otimes H^{-2n+2-m})$ as
$$
	D =
	    -  L_{e_a}
	      i(\overline{\vfld_{ab}})e^b,
 $$
where $L$ is the Lie derivative and $i$ is the interior multiplication.
The action of the horizontal form $e^b$ is the exterior multiplication.
\begin{lem}\label{D}
\begin{enumerate}
\item
The differential operator $D$ is invariant under the conformal automorphism
of $\R^{2n}$ if it is considered as an operator on $Z^\pm(\R^{2n})$.
\item\label{Dpend}
Let $l$ and $l'$ be non-negative integers. Then, we have
$$
	D \Gamma(\Zpm(M), \Lambda^{l,l'} \Zpm(M)\otimes H^{-2n+2-m})
		\subset \Gamma(\Zpm(M),
		  \Lambda^{l,l'} \Zpm(M)\otimes H^{-2n+2-m})
 $$
\item
We have $D^{n+1} = 0$.
\end{enumerate}
\end{lem}
\begin{pf}
We have (1) immediately by the definition.
Since $e^b$ and $L_{e_a}i(\overline{\vfld_{ab}})$ are commutative
and the space of horizontal $(0,1)$-forms are $n$-dimensional,
(3) follows by (2). We will prove (2) later in this section,
since it is needed a system of local coordinates of the fiber.
\end{pf}
 Put
$$
	F(x) = \sum_{k=0}^{\infty}\frac{x^k}{(k!)^2}.
 $$
An essential property of this function is the following lemma.
\begin{lem}\label{diff_eq}
Let $l$ be a non-negative integer. Then we have an equation
\begin{equation*}
	xF^{(l+2)}(x) =  - (l+1) F^{(l+1)}(x) + F^{(l)}(x).
\end{equation*}
\end{lem}

\begin{rem}
 The function $F$ also appears in the
construction of the inverse Penrose
transform of the Dolbeault complex ([I2]).
\end{rem}

Let $\Lambda^{0,\frac{1}{2}n(n-1)}_V$ be the line subbundle of
$\Lambda^{0,\frac{1}{2}n(n-1)}\Zpm(M)$ spanned by vertical forms. If we
identify $H^{-1}$ with $\overline{H}$ by the Hermitian
metric, we have
$$
	\Lambda^{0,\frac{1}{2}n(n-1)}\Zpm(M)\otimes H^{-2n+2-m} \supset
	  \Lambda^{0,\frac{1}{2}n(n-1)}_V\otimes H^{-2n+2-m} \simeq
	  \overline{H}^m,
 $$
since $\Lambda_V^{0,\frac{1}{2}n(n-1)}\simeq {\overline H}^{-2n+2}$.
By the theorem of Bott-Borel-Weil-Kostant,
we have an isomorphism
$$
	H^0(\Zpm, H^m)\simeq (S^m\D^{\pm})^*.
 $$
Hence we can define a lifting
$$
	j : \Gamma(M,S^m\D^{\pm}(M)) \rightarrow
	    \Gamma(\Zpm(M),
		\Lambda^{0,\frac{1}{2}n(n-1)}\Zpm(M)\otimes H^{-2n+2-m}).
 $$
By using $D$ and $j$,
we can define the local inverse Penrose transform as follows.
\begin{defn}\label{InvPenTr}
We define a map
\begin{align*}
	\AAA: \Gamma(M,S^m\D^{\pm}(M)) &\rightarrow
  	  \Gamma(\Zpm(M),
		\Lambda^{0,\frac{1}{2}n(n-1)}\Zpm(M)\otimes H^{-2n+2-m})\\
	\phi &\mapsto (n+m-2)! F^{(n+m-2)}(D)j(\phi).
\end{align*}
\end{defn}

\begin{rem}
We shall prove in the following section that
$\AAA$ does not depend on the conformally flat coordinates of $M$
and the trivialization of $V$
which are chosen. Hence $\AAA$ can be considered as a global
conformally invariant operator.
\end{rem}

Let us describe the Penrose transform
$\PPP_m$ by using the Dolbeault representation
of cohomology classes.
Since the complex dimension of $Z^\pm_x$ is
$\frac{1}{2}n(n-1)$, the isomorphism
$$
	H^{\frac{1}{2}n(n-1)}(\Zpm_x, H^{-2n+2-m})
	\simeq S^m\D^\pm_{\,x}
 $$
induces a map
\begin{equation*}
        \Gamma(\Zpm_x, \Lambda^{0,\frac{1}{2}n(n-1)}
             \Zpm_x\otimes H^{-2n+2-m})
	\rightarrow S^m\D^\pm_{\,x}.
\end{equation*}
Hence we have a map
$$
	\tilde{\PPP_m} :
	\Gamma(\Zpm(M), \Lambda^{0,\frac{1}{2}n(n-1)}
	     \Zpm(M)\otimes H^{-2n+2-m})
	\rightarrow
	\Gamma(M,S^m\Dpm(M)).
 $$
Then $\tilde{\PPP_m}$ induces the Penrose transform $\PPP_m$.
Since we have $\tilde{\PPP_m}(\AAA(\phi)) = \phi$ for any section
$\phi\in\Gamma(M,S^m\Dpm(M))$,
 we finish constructing the local inverse Penrose
transform by the following theorem.

\begin{thm}\label{CFlat}
Let $\phi$ be a section of $S^m\D^{\pm}(M)$.
Then $\AAA(\phi)$ is $\overline{\partial}$-closed, if
$d_m\phi=0$.
\end{thm}

The remainder of this section is devoted to a proof of this
theorem.
We prove it
by computing
 $dF^{(n+m-2)}(D)-F^{(n+m-2)}(D)d$.
\begin{lem}\label{commutator}
\begin{enumerate}
\item\label{exE}
Put $E = dD - Dd$. Then, we have
$$
	E = - L_{e_a} L_{\overline{\vfld_{ab}}}e^b.
 $$
\item\label{Gamma}
Let $d_H = e^aL_{e_a}$ be the exterior
differentiation to the horizontal direction, and put
$$
\Gamma = e^a\wedge e^bi(\overline{\vfld_{ab}})
		\sum_c(L_{\partdelx{c}})^2.
 $$
Then, we have
\begin{gather*}
	ED - DE = -2 Dd_H + \Gamma,\\
	\Gamma D - D\Gamma = 0,\\
	d_H D - D d_H = 0.
\end{gather*}
\item\label{commutatorItem}
Let $f(x)$ be a power series. Then, we have
$$
	df(D)-f(D)d = f'(D)E - f''(D)Dd_H + \frac{1}{2}f''(D)\Gamma.
 $$
\end{enumerate}
\end{lem}
\begin{pf}
(1)
Let $\Omega$ be the curvature form of $H^{-2n+2-m}$. Let $v$
be a vector field on $\Zpm(M)$. Then we have
$$
	[d, L_v] = - [i(v), \Omega]
 $$
as an operator acting on 
$\Gamma(\Zpm(M), \Lambda^{\!*}\Zpm(M)\otimes H^{-2n+2-m})$.
Since $[i(e_a),\Omega] = 0$ for any integer~$a$
such that $1\le a\le 2n$, we have the desired equation.\newline
\noindent (2) 
The second equation and the third one are
 immediate by the definitions and
the formulas
\begin{gather*}
	[L_v,i(v')] = i([v,v']),\\
	[L_v,L_{v'}] = L_{[v,v']} + \Omega(v,v'),
\end{gather*}
where $v$ and $v'$ are vector fields on $\Zpm(M)$.
By \eqref{lie_alg_str}, we compute
$$
	[\vfld_{ab},\vfld_{cd}] =
		- \delta_{ac}\vfld_{bd} + \delta_{ad}\vfld_{bc}
		+ \delta_{bc}\vfld_{ad} - \delta_{bd}\vfld_{ac},
 $$
so we have
$$
	[E,D] = -2Dd_H + \Gamma.
 $$
\noindent(3)
      By induction on $k$, we have
$$
	dD^k - D^kd = k D^{k-1}E - k(k-1) D^{k-1}d_H +
			\frac{1}{2}k(k-1)D^{k-2}\Gamma,
 $$
      which completes the proof.
\end{pf}

By Lemma \ref{diff_eq} and Lemma \ref{commutator} (\ref{commutatorItem}),
we have
\begin{multline*}
	dF^{(n+m-2)}(D) =\\ F^{(n+m-2)}(D)(d - d_H)
			   + F^{(n+m-1)}(D)(E + (n+m-1) d_H)
			+\frac{1}{2}F^{(n+m)}(D)\Gamma.
\end{multline*}
Since $j(\phi)$ is harmonic in the vertical direction, we have
$$
	F^{(n+m-2)}(D)(d-d_H)j(\phi) = 0.
 $$
If $\phi$ satisfy  $d_m\phi=0$, by Lemma \ref{commutator} (\ref{Gamma}),
we have
$$
	\frac{1}{2}F^{(n+m)}(D)\Gamma j(\phi) = 0.
 $$
Hence we complete the proof by computing the action of $E$.

Let us extend notation of a multi-index of the spin module in [I1],
which significantly simplifies computation as we will see below.
Let $(\theta_I)_{I < (1,\ldots,n)}$ be the basis of the spin module
$\D$ defined in [I1],
where $I<(1,\ldots,n)$ means that $I$ is a subsequence of
the sequence $(1,\ldots,n)$.
We regard a multi-index $I$ as a finite sequence of possibly
duplicate elements of $\{1,\ldots,2n\}$,
 and for $I=(i_1,\ldots,i_k)$,
$\theta_I$ is defined as
$	\theta_I = e_{i_1}*\cdots* e_{i_k}*\theta_{\emptyset}$
where $*$ is Clifford multiplication.
Let $(Z^I)_{I<(1,\ldots,n)}$ be the
 dual basis of $(\theta_I)$. Then we can consider $Z^I$ for
a multi-index $I$ in the same way.
Reduction of a multi-index is performed as follows.
Let $I$ be a reduced multi-index, i.e. $I$ is a subsequence of $(1,\ldots,n)$,
and let $i$, $j$ be distinct integers such that $1\le i,j\le n$.
Then
\begin{align*}
	\theta_{iiI} &= - \theta_I,\\
	\theta_{ijI} &= - \theta_{jiI},\\
	\theta_{(n+i)I} &=
	   \begin{cases}
	      \im \theta_{iI},&\text{if $i\not\in I$,}\\
	      -\im \theta_{iI},&\text{if $i\in I$.}
	   \end{cases}
\end{align*}
In the dual representation, we have
\begin{align*}
	Z^{iiI} &= - Z^I,\\
	Z^{ijI} &= - Z^{jiI},\\
	Z^{(n+i)I} &=
	   \begin{cases}
	     -\im Z^{iI},&\text{if $i\not\in I$,}\\
	      \im Z^{iI},&\text{if $i\in I$.}
	   \end{cases}
\end{align*}
Then we can reduce any multi-index $I$ to a unique reduced index
 $I'<(1,\ldots,n)$
such that $Z^I = \pm Z^{I'}$ or $Z^I = \pm\im Z^{I'}$ holds.
When a multi-index $I$ is used as a set (for example $i\in I$ and $I\cup J$),
it is considered
to be the set of numbers contained in the reduced form of $I$.
Let $|I|$ be the length of the reduced form of $I$.
Then we have
\begin{align*}
	\Dp &= \langle \theta_I \mid |I| \equiv 0 (2) \rangle,\\
	\Dm &= \langle \theta_I \mid |I| \equiv 1 (2) \rangle.
\end{align*}

Since $Z^\pm$ is a subvariety of $\Proj(\Dpm)$,
$Z^I$ can be considered as a homogeneous coordinate of $Z^\pm$.
For simplicity, we regard
$Z^I$ as a zero function when
 $|I|$  has
the parity opposite to that of projectivized spinors in which
the variety lies.
The defining equations of
$Z^\pm$ are given in [I1].
They can be given in our notation as follows.
\begin{lem}\label{spin_rel}
For multi-indices $I$ and $J$, let
 $d(I,J)=(I\setminus J) \cup (J\setminus I)$.
Let $a$ and $b$ be integers such that $1\le a,b\le 2n$.
Then we have the following relations on $Z^\pm$:
\begin{enumerate}
\item\label{spin_rel1}
$\displaystyle\sum_{k\in d(I,J)}Z^{kI}Z^{kJ} = 0.$
\item\label{spin_rel2}
$Z^{aI}Z^{aJ} = 0.$
\vspace{1mm}
\item\label{spin_rel3}
$\displaystyle\sum_{k\in d(I,J)}Z^{akI}Z^{bkJ} =
		Z^IZ^{abJ} - Z^{abI}Z^J$.
\end{enumerate}
\end{lem}
\begin{pf}
The equation (1) is an immediate consequence of [I1] Corollary 3.3.
Let $i$ be an integer such that $1\le i \le n$. Then
$$
	Z^{iI}Z^{iJ} + Z^{(n+i)I}Z^{(n+i)J} =
		\begin{cases}
			2Z^{iI}Z^{iJ}, & i\in d(I,J),\\
			0, & i\not\in d(I,J).\\
		\end{cases}
 $$
Hence (2) follows immediately by (1).
We can prove (3) by simple computation using~(2).
\end{pf}
Now we fix a multi-index $I$.
Let $z^J = Z^J/Z^I$ for a multi-index $J$.
We let $w_{ij}=z^{ijI}$. Then $(w_{ij})_{1\le i<j\le n}$ are local
coordinates on
$U_I=\{[Z^J]\in \Zp \cup \Zm \mid Z^I\ne 0\}$.

\begin{lem}\label{localdeliv}
For integers $i$ and $j$ such that $1\le i<j\le n$, we have
$$
	\partdel{z^J}{w_{ij}} =
	  \begin{cases}
	    z^{jiJ},& i,j\in d(I,J),\\
 	    0,&\text{otherwise}.
	  \end{cases}
 $$
\end{lem}
\begin{pf}
We have the relation
\begin{equation}\label{four_rel}
	z^J = (z^{abJ}z^{cdJ}-z^{acJ}z^{bdJ}+z^{adJ}z^{bcJ})/z^{abcdJ},
\end{equation}
by Lemma \ref{spin_rel} \eqref{spin_rel3}.
Thus we can prove the lemma
inductively
by taking appropriate integers $a$, $b$, $c$, and $d$.
\end{pf}
\begin{lem}\label{vectorfield}
The vector field $\vfld_{ab}$ is written
in the local coordinates as
$$
	\vfld_{ab} = \sum_{1\le i<j\le n} \frac{1}{2} 
			(z^{bjI}z^{aiI} - z^{ajI}z^{biI})\partdel{}{w_{ij}}.
 $$
\end{lem}
\begin{pf}
Since the one-parameter subgroup of $\Spin(2n)$
corresponding to $F_{ab}$ is $
\cos \frac{t}{2} + \sin \frac{t}{2}\, e_ae_b$, 
its action is written as
$$
	Z^J(t) = \cos \frac{t}{2} Z^J + \sin \frac{t}{2} Z^{baJ}.
 $$
Thus,  by using \eqref{four_rel}, we compute
$$
	\left. \frac{dw_{ij}(t)}{dt} \right|_{t=0}
	= \frac{1}{2}(z^{bjI}z^{aiI} - z^{biI}z^{ajI}),
 $$
which completes the proof.
\end{pf}

\begin{lem}\label{horizontal_form}
The space of horizontal $(1,0)$-forms on $\Zpm(M)$ is spanned by
$$
	\alpha^J = -z^{aJ}e^a, \qquad J<(1,\ldots,n).
 $$
\end{lem}

\begin{pf}
It suffices to show the lemma when $M=\R^{2n}$.
Its Riemannian twistor space is given in [I1] \S 8.
Let $\D'$ be a spin module of $\Spin(2n+2)$:
\begin{equation*}
	\D' = \langle \theta_J \mid J< (0,1,\ldots,n) \rangle.
\end{equation*}
Let ${Z'}^+\subset\Proj({\D\!'}^+)$ be the
parameter space of the compatible complex structure
of the vector space $\R^{2n+2}$,
 which can be identified with the twistor space $Z^+(S^{2n})$.
Since  stereographic projection defines a conformal embedding
$\R^{2n} \subset S^{2n}$, $Z^+(\R^{2n})$
 is an open submanifold of $Z^+(S^{2n})$.
Let $(Z^J)$ be the homogeneous coordinates with respect to $(\theta_J)$.
Then we have:
\begin{equation*}
	Z^+(\R^{2n}) = \{ (Z^J)_{J<(0,\ldots,n)} \in {Z'}^+ \mid
			\exists J < (1,\ldots,n) \quad\hbox{such that} 
			\quad Z^J \not= 0 \}.
\end{equation*}
Since a translation of $\R^{2n}$ is a conformal transformation,
 it induces a
holomorphic transformation of $Z^+(\R^{2n})$,
 which is representable by an element of $\Spin(2n+2;\C)$.
Let $x = (x_1,\ldots,x_{2n})$ be an element of $\R^{2n}$.
Then the corresponding element of $\Spin(2n+2;\C)$ is:
$$
	\alpha(x) = 1 + \frac{1}{2} \sum_{a=1}^{2n}
			 x_a e_a (\im e_0 + e_{0'}),
 $$
where we think the standard basis of $\R^{2n+2}$ to be
$(e_0,e_1,\ldots,e_n,e_{0'},e_{n+1},\ldots,e_{2n})$.
A point on the fiber over $0\in\R^{2n}$ is written as
$\sum_{I\not\ni 0} Z^I\theta_I$,
thus its image by the transformation $\alpha(x)$ is written as
$$
	\alpha(x)\sum_{J\not\ni0}Z^J\theta_J =
		Z^J\theta_J +
		\im \sum_{a=1}^{2n} x_a Z^{aJ}\theta_{0J}.
 $$
Hence, for a multi-index $J$ such that $0\not\in J$,
 the homogeneous function
\begin{equation*}
	Z^{0J} = \im \sum_{a=1}^{2n} x_a Z^{aJ}
\end{equation*}
is holomorphic, which completes the proof on $Z^+(M)$.
The proof on $Z^-(M)$ is done in the same way.
\end{pf}

\begin{rem}
By the definition of the almost complex structure of twistor
spaces,
the lemma is true for any spin manifold $M$ and any
oriented orthonormal
frame $(e^a)$ of~$T^*M$.
\end{rem}

Now we can calculate the action of $E$.
Let $I_1,\ldots,I_m$ be multi-indices having the same parity of length.
Let $\symtheta\in S^m\Dpm$ be the symmetrization of
$\theta_{I_1}\otimes\cdots\otimes\theta_{I_m}$.
Then we define
$$
	\sss = j(\symtheta).
 $$
We compute its action on $M\times U_I$ for a multi-index $I$.
If $m\ge1$, we assume $|I| \equiv |I_1| (2)$, since
$U_I = \emptyset$ in other cases.

\begin{lem}\label{comp_E}
We have
\begin{multline*}
	\left(E+ (n+m-1) d_H\right)
	 j(\phi) =\\
		 -\frac{m}{2} \partdel{
			\phi^{abI_1,I_2,\ldots,I_m}}{x_a} e^b\wedge\sss
		  - \frac{2n-2+m}{2N}
			\partdel{\phi^{\symindex}}{x_a}
			\overline{z^{aJ}}\alpha^J\wedge
			\sss.
\end{multline*}
where $N = \sum_{J\not\ni0} |z^J|^2$ is the Hermitian metric of $H^{-1}$
with the standard trivialization on $M\times U_I$.
\end{lem}
\begin{pf}
Let $a$ and $b$ be distinct
  integers such that $1\le a, b\le 2n$. Then we claim
\begin{equation}\label{dirt}
	L_{\overline{\vfld_{ab}}}\sss = \frac{1}{2}
	  (- \frac{2n-2+m}{N}\sum_{J}z^{bJ}\overline{z^{aJ}}
		 + \sum_{i=1}^{m}\frac{\overline{z^{baI_i}}}
				     {\overline{z^{I_i}}})
	\sss,
\end{equation}

Let $\overline{\rho^J}$ be the section
of $\overline{H}$ corresponding to $\theta_J$.
Let $\overline{K_I}$ be the standard trivialization of 
$\Lambda^{0,\frac{1}{2}n(n-1)}_V$
over $M\times U_I$. Then, by the definition of $j$, we can write
\begin{equation*}
	\overline{s^{\symindex}} = \overline{K_I}\otimes
		{\overline{\rho^I}}^{2n-2}\otimes
		\overline{\rho^{I_1}}\otimes\cdots\otimes
		\overline{\rho^{I_m}}.
\end{equation*}
First, we compute
\begin{align*}
	L_{\overline{\vfld_{ab}}} \overline{\rho^I} &= 
		\nabla_{\overline{\vfld_{ab}}} \overline{\rho^I}\\
		&= -\frac{\overline{\vfld_{ab}}(N)}{N} \overline{\rho^I}\\
		&= -\sum_{1\le i\not=j\le n}
			\frac{\overline{z^{bjI}z^{aiI}}}{2N}
				(\sum_{J\text{ s.t. } i,j\in d(I,J)}
				  z^J\overline{z^{jiJ}})
			\overline{\rho^I}
			\qquad\text{[By Lemma \ref{localdeliv}]}\\
		&= -\frac{1}{2N}
			\sum_J\sum_{j\in d(I,J)}
				\overline{z^{bjI}z^{ajJ}}z^J
			\overline{\rho^I}
			\qquad\text{[By Lemma \ref{spin_rel}
				\eqref{spin_rel3}]}\\
		&= \frac{1}{2}(-\frac{1}{N}\sum_Jz^{bJ}\overline{z^{aJ}}
		   + \overline{z^{baI}})
			\overline{\rho^I}
			\qquad\text{[By Lemma \ref{spin_rel}
				\eqref{spin_rel3}]}.
\end{align*}
Hence, for another multi-index $I'$ such that
$|I'|\equiv |I| (2)$,  we have
$$
	L_{\overline{\vfld_{ab}}} \overline{\rho^{I'}} = 
		\frac{1}{2}(-\frac{1}{N}\sum_Jz^{bJ}\overline{z^{aJ}}
		   + \frac{\overline{z^{baI'}}}{\overline{z^{I'}}})
			\overline{\rho^{I'}}.
 $$
Second, we compute
$$
	L_{\overline{\vfld_{ab}}}\overline{K_I} =
		- (n-1)\overline{z^{baI}K_I}.
 $$
Therefore we have the equation \eqref{dirt}.

 By Lemma \ref{commutator} \eqref{exE} and \eqref{dirt}, we compute
\begin{align*}
	E j(\phi) 
	 &= - \sum_{a\not=b} \frac{1}{2}
		\partdel{\phi^{\symindex}}{x_a}e^b
		\left(-\frac{2n-2+m}{N}\sum_J z^{bJ}\overline{z^{aJ}}
                    + \sum_{i=1}^{m}
                        \frac{\overline{z^{baI_i}}}{\overline{z^{I_i}}}
                        \right)
			\overline{s^{\symindex}}\\
	 &=  -\frac{m}{2}
		  \partdel{\phi^{abI_1,I_2,\ldots,I_m}}{x_a}
		    e^b\wedge\overline{s^{\symindex}}
		  - \frac{2n-2+m}{2N}
			\partdel{\phi^{\symindex}}{x_a}
			\overline{z^{aJ}}\alpha^J\wedge
			 \sss\\
	 &\phantom{=}\,
	    - (n-1+m)d_H j(\phi).
\end{align*}
Hence we complete the proof.
\end{pf}

\begin{pf*}{\it Proof of Lemma \ref{D} \eqref{Dpend}}
We have
$$
	D 
	  = \frac{1}{2}\sum_{1\le i<j\le n}
		\idwij (\overline{z^{aiI}\alpha^{jI}}
			- \overline{z^{ajI}\alpha^{iI}})
		L_{e_a}.
 $$
Since $L_{e_a}(dw_{ij}) = 0$,  $L_{e_a}(\alpha^J) = 0$ and
$L_{e_a}(\overline{\rho^J})=0$,
we complete the proof.
\end{pf*}

Since $D$ and $\alpha^J$ are commutative,
 the second term of Lemma \ref{comp_E}
can be neglected  modulo $(1,0)$-forms.
Hence, if $m=0$, we  finish the proof of the theorem.
If $m>0$, the coefficient of the first term 
is that of $\theta_{bI_1}\otimes\theta_{I_2,\ldots,I_m}$
in $\frac{m}{2}d_m(\phi)$.
Hence, if $d_m(\phi)$ vanishes, we have
$$
	(E + (n+m-1)d_H)j(\phi) \equiv 0\qquad \text{modulo $(1,0)$-forms,}
 $$
for any non-negative integer $m$.
Thus we complete the proof of the theorem.

\section{Well-definedness of $\AAA$ as a global operator}

By Theorem \ref{CFlat}, we prove that
Definition \ref{InvPenTr} gives an inverse
Penrose transform 
locally when the base
manifold is  conformally flat. In this section we show
that the construction is
independent with respect to the conformally flat coordinates
which are chosen and
gives the global inverse Penrose transform.

We continue assuming that the base manifold
is $2n$-dimensional with $n\ge 3$,
hence the metric of the base manifold is conformally flat.
We have a local inverse Penrose transform $\AAA$ by
Definition \ref{InvPenTr} on each chart which has conformally
flat coordinates.
By a theorem of Liouville ([D.F.N] Theorem 15.2),
a coordinate transformation is a
certain restriction of an orientation preserving
conformal automorphism of $S^{2n}$ fixing the point which
is regarded as the center of each chart.
The orientation preserving
conformal automorphism group of $S^{2n}$ is $\SO_0(1,2n+1)$.
Let $G$ be the isotropy subgroup at $0\in\R^{2n}\subset S^{2n}$.
Then we shall show the invariance of  $\AAA$
 under the action of $G$ near the origin.

Let $\CE(2n)$ be the orientation preserving
conformal automorphism group  of $\R^{2n}$.
Since  a local conformal map
can be  extended uniquely to a conformal automorphism
of $S^{2n}$, $\CE(2n)$ can be considered to be
 a subgroup of $\SO_0(1,2n+1)$.
Let $\CO(2n)$ be the isotropy subgroup of $\CE(2n)$ at the origin.
Let $\tau$ be
the conformal map defined as
\begin{align*}
	\tau: \R^{2n}\setminus \{0\} &\rightarrow \R^{2n}\setminus \{0\}\\
	x &\mapsto \frac{x}{|x|^2}.
\end{align*}
By computing the Lie algebra of $G$, we can show that
the group $G$ is generated by $\CO(2n)$ and
$\tau \circ T(x) \circ \tau$ for
$x\in\R^{2n}$, where $T(x)$ is
the translation map on $\R^{2n}$ by~$x$.
Since $D$ and $j$ are invariant under the action of $\CE(2n)$,
 $\AAA$ is invariant under the action of $\CO(2n)$.
It is also invariant under $T(x)$
for any $x\in\R^{2n}$,
 so we can show its invariance under
 $\tau\circ T(x)\circ\tau$
on $\R^{2n}\setminus  \{0,\tau(-x)\}$
by showing its invariance under~$\tau$ on $\R^{2n}\setminus\{0\}$.
This proves the invariance of $\AAA$
under $\tau\circ T(x)\circ\tau$ near the origin,
since $\AAA(\phi)$ is expressible
as  a polynomial of jets of $\phi$.

Hence it suffices to show the invariance of $\AAA$ under~$\tau$
 on $\R^{2n}\setminus \{0\}$.
The following lemma is used to reduce the computation
to  a certain point on the fiber over the point
$x_0 = {}^t(1,0,\ldots,0)$ instead of computing it on the whole space
 $Z^\pm(\R^{2n}\setminus \{0\})$.
\begin{lem}\label{trans}
The group $\CO(2n)$ acts transitively on $Z^{\pm}(\R^{2n}\setminus\{0\})$.
\end{lem}
\begin{pf}
Since $\CO(2n)$ acts transitively on $\R^{2n}\setminus\{0\}$,
 it suffices to prove the transitivity on the fiber
over the point $x_0$.
The isotropy subgroup at $x_0 $
contains the subgroup which is naturally identified with $\SO(2n-1)$.
Then it acts on $Z^\pm_{x_0}$,
and the isotropy subgroup
at $(x_0,[\theta_I])$
for an appropriate $I$ is naturally identified with
$\Un(n-1)$. Thus the isomorphism
$$
	\SO(2n-1)/\Un(n-1) \simeq \SO(2n)/\Un(n)
 $$
means that $\SO(2n-1)$ acts  transitively on the fiber
$Z^\pm_{x_0}$.
\end{pf}

Let $B$ be an endomorphism defined as
$$
	B = - [D,|x|^2],
 $$
where $|x|^2$ is considered to be an operator
by the multiplication.
\begin{lem}\label{DxB}
The operators $D$ and $B$ are commutative.
\end{lem}
\begin{pf}
We fix a multi-index $I$.
Let $z^J$ and $w_{ij}$ be
as in the previous section. Then, by Lemma \ref{vectorfield}
and Lemma \ref{spin_rel} \eqref{spin_rel2}, we have
\begin{align*}
	[D,B] 
	      &= - \frac{1}{2}
		\sum_{a\not=b,d}
		\sum_{i\not=j}
		\sum_{k\not=l}
		\overline{z^{bjI}z^{aiI}z^{dlI}z^{akI}}
		\idwij e^b \idwkl e^d\\
	      &= \frac{1}{2}
		\sum_b
		\sum_d
                \sum_{i\not=j}
		\sum_{k\not=l}
		\overline{z^{bjI}z^{dlI}(z^{biI}z^{bkI}+z^{diI}z^{dkI})}
		\idwij e^b \idwkl e^d,
\end{align*}
where we put $\partdel{}{w_{ji}} = - \partdel{}{w_{ij}}$
for integers $i$ and $j$ such that $1\le i < j\le n$.
Since, for fixed integers $b$ and $d$, we have
$$
	\sum_{i\not=j}
	 \overline{z^{bjI}z^{biI}}\idwij = 0,\quad
	\sum_{k\not=l}
 	 \overline{z^{dlI}z^{dkI}}\idwkl = 0,
 $$
we complete the proof.
\end{pf}
We can relate
$\tau^*(j(\phi))$ and $j(\tau^*\phi)$
by using $B$  as follows.
\begin{lem}\label{pullbackj}
We have $\tau^*(j(\phi)) = \exp(|x|^{-2}B) j(\tau^*\phi)$.
\end{lem}
\begin{pf}
Let $\kappa$ be the orientation reversing isometry defined as
$$
	\kappa\left({}^t(x_1,x_2,\ldots,x_{2n})\right) =
		{}^t(-x_{1},x_{2},\ldots,x_{2n}).
 $$
Let $G'$ be the transformation group generated by
$\kappa$ and $\CO(2n)$. Then we have
\begin{enumerate}
\item an element of $G'$ preserves $j$ and $\frac{1}{|x|^2}B$,
\item $\tau G' = G'\tau$,
\item $G'$ acts on
  $\Zp(\R^{2n}\setminus\{0\}) \cup \Zm(\R^{2n}\setminus\{0\})$
transitively.
\end{enumerate}
Hence it suffices to compute them at the point
$z_0 = (x_0,[\theta_{\emptyset}])$.
The map between the twistor spaces
induced by $\tau$ is written in the homogeneous coordinates as
\begin{equation*}
\begin{split}
	\tau : Z^+(\R^{2n}\setminus\{0\}) & \rightarrow
		Z^-(\R^{2n}\setminus\{0\})\\
	(x,[Z^I]) & \mapsto
		\left(\frac{x}{|x|^2},[x_aZ^{aI}] \right).
\end{split}
\end{equation*}
Hence $\tau(x_0,[\theta_{\emptyset}]) =
(x_0,[\theta_{1}])$. We take the two systems of
local coordinates
\begin{equation*}
\begin{matrix}
	w_{ij}  &= {\displaystyle\frac{Z^{ij}}{Z^\emptyset}},
		\qquad &\text{on $U_\emptyset$,}\\[10pt]
	w'_{ij} &= {\displaystyle\frac{Z^{ij1}}{Z^1}},
		\qquad & \text{on $U_1$,}
\end{matrix}
\qquad 1\le i<j \le n.
\end{equation*}
Then we have
$$
	\left. \tau^*\overline{dw'_{ij}}\right|_{z_0} = 
	\begin{cases}
		-\overline{dw_{1j}} + e^j - \im e^{n+j},& i=1,\\
		\overline{dw_{ij}},& i>1.
	\end{cases}
 $$
For each $j$, we have
$$
	\left(
		(e^j - \im e^{n+j})
			\idw{1j}
	\right)^2 = 0.
 $$
Thus we have
$$
	\prod_{j=2}^{n} \left(
		1 - (e^j - \im e^{n+j})
			\idw{1j}
			\right)
	= \left.\exp\left(|x|^{-2}B\right)\right|_{z_0}.
 $$
Hence
$\left.\tau^{*}\Lambda_V^{0,\frac{1}{2}n(n-1)}Z^{-}(\R^{2n}\setminus\{0\})
\right|_{z_0}$
is spanned by
$$
	\exp\left(|x|^{-2}B\right)
		\overline{dw_{12}}\wedge\cdots
			\wedge\overline{dw_{n-1,n}}.
 $$
Since we have
 $\tau^{*}(j(\phi)) \equiv j(\tau^{*}\phi)$
modulo horizontal forms,
we complete the proof.
\end{pf}

\begin{lem}\label{Vanish}
\begin{enumerate}
\item\label{vs}
We have $D\overline{s^{I_1,\ldots,I_m}} = 0$.
\item\label{vt}
Put $\tttt=|x|^{2(n+m-1)}j(\tau^{*}\theta_{\symindex})$.
Then we have $D\overline{t^{I_1,\ldots,I_m}} = 0$.
\end{enumerate}
\end{lem}
\begin{pf}
Since $L_{e_a}\sss = 0$ for an integer $a$ such that
$1\le a\le 2n$, we have (\ref{vs}).

Let $J$ be a multi-index.  By using Lemma \ref{spin_rel} \eqref{spin_rel2},
 we have
$$
	[D, x_{a}\overline{z^{aJ}}] 
	   = \frac{1}{2}\sum_b \sum_{i\not=j}\overline{z^{bjI}z^{biI}z^{bJ}}
		\idwij e^b
	   = 0.
 $$
The $m$-th symmetric spin bundle $S^m\D^{\pm}(\R^{2n}\setminus \{0\})$
with conformal weight $n-1+\frac{m}{2}$
is transformed by $\tau$
as follows.
\newcommand{\tmpsign}{c}
\begin{equation*}
	\tau^{*}\theta_{\symindex} =
		\tmpsign
		 |x|^{-2(n+m-1)}x_{a_1}\cdots x_{a_m}
			\theta_{a_1I_1,\ldots,a_mI_m},
\end{equation*}
where $c$ is the constant number determined by
 spin structures of both ends of $\tau$.
Then
$$
	\tttt 
	 = \tmpsign (x_{a_1}\overline{z^{a_1I_1}})\cdots
		    (x_{a_m}\overline{z^{a_mI_m}})
		\overline{s^{\emptyset,\ldots,\emptyset}}.
 $$
Hence we complete the proof of (\ref{vt}).
\end{pf}

We have $\tau \CO(2n) = \CO(2n)\tau$. 
We have also that elements of $\CO(2n)$ preserve $\AAA$.
Hence, by Lemma \ref{trans},
it suffices to show that
$\AAA$
is invariant under $\tau$ at a certain point on the fiber over $x_0$.
Actually, we do not need to specify a special point on the fiber,
so  we compare two inverse Penrose transforms
 on the fiber over the point $x_0$.
For simplicity, we write $x'=(x_{2},\ldots,x_{2n})$.
Since they are linear with respect to $\phi$,
it suffices to show $\tau^{*}(\AAA(\phi)) = \AAA(\tau^{*}\phi)$
 with respect to
the section
$$
	\phi= {x_{1}}^l f(x')\theta_{I_1,\ldots,I_m}
 $$
where $l$ is a non-negative
integer and $f(x')$ is a homogeneous polynomial of degree~$l'$.

\begin{lem}\label{pullback}
\begin{enumerate}
\item\label{DxXzero}
Let $a$ be an integer such that
$1\le a\le 2n$. Then
 $[D,x_a]$ and $D$ are commutative.
\item\label{DxXone}
At points of the fiber over $x_0$,
$D^{l'}f(x')$ is an endomorphism of the vector bundle,
 and the pull-back by $\tau$ is computed as
$$
	\tau^{*}\left(\left.D^{l'}f(x')\right|_{p^{-1}(x_0)}\right)
		 = D^{l'}f(x').
 $$
\item\label{DxXtwo}
The pull-back of $B$ by $\tau$ is computed as
$$
	\left.\tau^{*}B\right|_{p^{-1}(x_0)} = - B.
 $$
\end{enumerate}
\end{lem}
\begin{pf} We have
\begin{equation}\label{Dxa}
	[D, x_{a}] = - i(\vfld_{ab})e^b.
\end{equation}
Since this is constant with respect to $x_1,\ldots,x_{2n}$, we
prove the assertion~(1).
The Jacobian matrix of
$\tau : \R^{2n}\setminus \{0\} \rightarrow \R^{2n}\setminus \{0\}$  is
$$
	|x|^{-2}R(x),
 $$
where $R(x)$ is the reflection with respect to the hyperplane
with the normal vector~$x$.
Since the operator (\ref{Dxa})
is transformed as a one-form, we have
\begin{equation*}
	\left.\tau^{*}[D, x_{a}]\right|_{p^{-1}(x_0)} =
		\begin{cases}
			-[D, x_{1}] & a = 1,\\
			 \,[D,x_{a}] & a > 1.
		\end{cases}
\end{equation*}
Hence we have the equation of \eqref{DxXone}.
Since $B|_{p^{-1}(x_0)}
 = - 2[D,x_{1}]$, we also have the equation of (\ref{DxXtwo}).
\end{pf}
If $k<l'$, then we have $D^k f(x')=0$
on the fiber over $x_0$. 
Hence, by Lemma \ref{pullback} \eqref{DxXone}, we have
$$
	\tau^{*} \left(
		\AAA(\phi)
		  \right)
	 = \frac{(n+m-2)!}{l'!}
		    D^{l'}f(x')
		    \tau^{*}
		      \left(
			F^{(n+m-2+l')}(D) {x_{1}}^l
			\overline{s^{I_1,\ldots,I_m}}
		      \right).
 $$
On the other hand, since ${x_{1}}^l f(x')$ is homogeneous
of degree $l+l'$, we have
$$
	\tau^{*}\phi =
		|x|^{-2(l+l')} {x_{1}}^l f(x')
		\tau^{*}\theta_{\symindex}.
 $$
In the same way, we have
$$
	\AAA(\tau^{*}\phi)
	=	\frac{(n+m-2)!}{l'!}
		D^{l'}f(x')
		F^{(n+m-2+l')}(D)
		  |x|^{-2(n+m-1+l+l')} {x_{1}}^l\tttt.
 $$
Hence we can show the invariance of $\AAA$ under $\tau$
by the following lemma.
\begin{lem} Let $n'$ be a non-negative integer. At  points
on the fiber over $x_0$, we have
$$
	\tau^{*}\left(
		F^{(n')}(D){x_{1}}^l
		\overline{s^{I_1,\ldots,I_m}}
	    \right) =
	 F^{(n')}(D) |x|^{-2(n'+1+l)} {x_{1}}^l
 		\tttt.
 $$	
\end{lem}
\begin{pf}
We prove this by induction on $l$.
Let $l=0$. 
By Lemma \ref{Vanish} (\ref{vs}) and Lemma~\ref{pullbackj}, we have
$$ 
	n'!\, \tau^{*}\left( F^{(n')}(D)\sss\right) 
	 = \exp(B)j(\tau^{*}\symtheta).
 $$
On the other hand, by Lemma \ref{DxB} and Lemma \ref{Vanish} (\ref{vt}),
we have
$$
	n'!\, F^{(n')}(D) |x|^{-2(n'+1)}\tttt
	 = \exp(B)j(\tau^{*}\symtheta).
 $$
Hence we have the equation when $l=0$.

Let us assume that the equation is satisfied for 
integers less than or equal to  $l$. Then
\begin{align*}
	F^{(n')}(D){x_{1}}^{l+1}\sss
	&= \sum_{k=0}^{\infty}
		\frac{1}{k! (n'+k)!}
		D^kx_1{x_{1}}^{l}\sss\\
	&= \sum_{k=0}^{\infty}
		\frac{1}{k! (n'+k)!}
		\left(
			-\frac{k}{2} BD^{k-1} + D^k
		\right)
		{x_{1}}^l\sss\\
	&= - \frac{1}{2} BF^{(n'+1)}(D){x_{1}}^l\sss
		+ F^{(n')}(D){x_{1}}^l\sss.
\end{align*}
Hence, by Lemma \ref{pullback} (\ref{DxXtwo})
and the hypothesis, we have
\newcommand{\xone}{{x_{1}}}
\newcommand{\pullbS}{x_{a_1}\cdots x_{a_m}
	\overline{s^{a_1I_1,\ldots,a_mI_m}}}
\begin{align*}
	&\tau^{*}\left(F^{(n')}(D){x_{1}}^{l+1}\sss\right)\\
	&\quad =
		\frac{1}{2} B
		F^{(n'+1)}(D)|x|^{-2(n'+2+l)}\xone^l\tttt
		+ F^{(n')}(D)|x|^{-2(n'+1+l)}\xone^l\tttt.
\end{align*}
Since
$$
	F^{(n')}(D) (\xone - |x|^2)|x|^{-2(n'+2+l)}\xone^l
		\tttt
	=
		\frac{1}{2} BF^{(n'+1)}(D)
		|x|^{-2(n'+2+l)}\xone^l \tttt,
 $$
we complete the proof.
\end{pf}

Hence we have the equation
$\AAA(\tau^{*}\phi) = \tau^{*}\left(\AAA(\phi)\right)$
for any section $\phi$.
Thus we have proved the remark after Definition \ref{InvPenTr}.

\section{The inverse Penrose transform over four-manifolds}

In this section, as an extension of Definition \ref{def4},
we give the inverse Penrose transform over
a four-manifold.
It is an interpretation of the formula given by
Woodhouse. 

Let $M$ be a four dimensional
spin manifold, and let $V$ be a Hermitian vector
bundle on $M$ with a connection.
Since reversing the orientation of $M$ exchanges $\Zp(M)$ and $\Zm(M)$,
it suffices to consider the inverse Penrose transform
 only on $Z^+(M)$.
Thus we assume that the Riemannian metric of $M$ is anti-self-dual.
Assume also that the connection of $V$ is anti-self-dual. 
This means that the pull-back $p^{-1}V$ can be
naturally considered to be a holomorphic vector bundle
with a holomorphic connection on the complex
manifold $Z^+(M)$.

Let us define a differential operator $D$.
In conformally flat case, we have local
conformally flat coordinates, which significantly simplify
the computation.
But we do not have such coordinates on
anti-self-dual manifolds.
So we shall define $D$ by using an arbitrary local orthonormal frame of $TM$.

Let $(e_a)$ be an oriented
 local orthonormal frame of $TM$
 on an open subset $U$. 
Let $(e^a)$ be the dual frame of $T^*M$.
Since the twistor space is a fiber bundle over $M$
associated to the orthonormal frame bundle,
we can define horizontal tangent vectors
by the Levi-Civita connection.
Therefore, we consider $(e_a)$ as
a local frame of the space of horizontal vector fields.
Also $(e^a)$ is regarded as a local frame
of the space of horizontal forms.
Let $\omega$ be the connection form of 
$TU$ with respect to the Levi-Civita connection.
For each integer $a$ such that $1\le a \le 2n$,
we define a differential operator acting on
$\Gamma(\Zp(U), \Lambda^{*}Z^+(U)\otimes H^{-2-m}\otimes p^{-1}V)$ by
\newcommand{\lhat}[1]{\hat{L_{e_{#1}}}}
\newcommand{\lhatvar}[1]{\hat{L_{{e'}_{#1}}}}
$$
	\lhat{a} = L_{e_a} + i(e_b)\omega^b_a.
 $$
Then, the following lemma can be proved by simple computations.
\begin{lem}\label{trL}
Let ${e'}_a = e_bh^b_a$ be another local orthonormal frame.
Then we have $\lhatvar{a}=\lhat{b}h^b_a$.
\end{lem}

With respect to the local trivialization $(e_a)$,
we define a $(1,0)$-vector field $\vfld_{ab}$ on $Z^+(U)$
by using the $\SO(4)$-action on $Z^+(U)$ as in \S\ref{const}. 
The next lemma is an immediate consequence of the definition.

\begin{lem}\label{trF}
Let $({e'}_a)$ be as above.
Let $\vfld'_{ab}$ be the vector filed
corresponding to the frame $({e'}_a)$.
Then we have
$\vfld'_{ab} = {h^{-1}}^a_c\vfld_{cd}h^d_b$.
\end{lem}

Let us define an operator $D$ acting on
$\Gamma(\Zp(U), \Lambda^{*}Z^+(U)\otimes H^{-2-m}\otimes p^{-1}V)$
as follows
$$
	D = -\lhat{a}i(\overline{\vfld_{ab}})e^b.
 $$
Then we have the following
lemma by Lemma \ref{trL} and Lemma \ref{trF}.
\begin{lem}
The operator $D$ is defined
independently of the choice of a local
orthonormal frame of $TU$.
Hence it is considered to be a global
 operator acting 
on $\Gamma(\Zp(M), \Lambda^{*}Z^+(M)\otimes H^{-2-m}\otimes p^{-1}V)$.
\end{lem}

The Levi-Civita connection on the base manifold
defines the decomposition of the cotangent bundle
of the twistor space.
Hence, as in \S\ref{const},
we define a map:
$$
        j : \Gamma(M,S^m\Dp(M)\otimes V) \rightarrow
            \Gamma(\Zp(M),
                \Lambda^{0,1}\Zp(M)\otimes H^{-2-m}
		\otimes p^{-1}V).
 $$
Since the decomposition is global in this case,
the map $j$ is also global.

\begin{defn}\label{def4}
A map $\AAA$ is defined as
\begin{align*}
	\AAA: \Gamma(M,S^m\Dp(M)\otimes V) &\rightarrow
		\Gamma(Z^+(M),
		  \Lambda^{0,1}\Zp(M)\otimes H^{-2-m}\otimes p^{-1}V)\\
		\phi &\mapsto j(\phi) + \frac{1}{m+1}Dj(\phi).
\end{align*}
\end{defn}

It can be shown by straightforward computation that 
this differential form is equivalent to
the inverse Penrose transform given by Woodhouse ([W], \S5),
when $V$ is a trivial line bundle.
Computations in [W] can be easily extended to
the case of non-trivial bundle $V$. Hence we have the following theorem.

\begin{thm}[{[W]}]\label{lastth}
Let $\phi$ be a solution of $d_m\phi=0$,
 then $\overline{\partial}\AAA(\phi)=0$.
Hence $\AAA$ gives the inverse Penrose transform.
\end{thm}

\noindent {\it Remarks.}
 1. The transform $\AAA$ does not depend on
the Riemannian structure but the conformal structure of $M$.
Hence the above definition is equal to
Definition \ref{def4}
in the case of a flat vector bundle over
a conformally flat four dimensional manifold.

2. Since we have $D^2j(\phi)=0$,
we can write
$$
	\AAA(\phi) = m! F^{(m)}(D)j(\phi),
 $$
as in  the higher dimensional cases.


\end{document}